\documentstyle[12pt,leqno]{article}

\newtheorem{Definition}{Definition}[section]
\newtheorem{Proposition}{Proposition}[section]

\def\textmap#1{\mathop{\vbox{\ialign{
                                ##\crcr
    ${\scriptstyle\hfil\;\;#1\;\;\hfil}$\crcr
    \noalign{\kern-1pt\nointerlineskip}
    \rightarrowfill\crcr}}\;}}

\date{December 14, 1995}

\begin{document}

\begin{flushright}
alg-geom/9512011
\end{flushright}
\vspace{1.4 cm}

\begin{center}
{\large{\bf On Invariant Theory }} \\
\vspace{2.2 cm}
 Valeri V. DOLOTIN \\
\vspace{0.8 cm}
{\it vd@main.mccme.rssi.ru} \\College of Mathematics, Independent
University of Moscow \\ \vspace{1.4 cm}

\begin{abstract}
Here  we develop a technique of computing the invariants of  $n-$ary forms
and systems of forms using the discriminants of  corresponding multilinear
forms buit of their partial derivatives, which  should  be cosidered  as
generalizations of classical discriminants and resultants for binary
forms.
\end{abstract}
\end{center}

\setcounter{section}{-1}

\section{Introduction}

Let $f(y_0,y_1,...,y_n)=\displaystyle
\sum^{}_{\alpha_1\le{...}\le{\alpha_d}}c_{\alpha_1{...}{\alpha_{d}}%
}y_{\alpha_{1}} {...}y_{\alpha_{d}}$ be an $n+1-$ary form of degree $d$,
which is a polynomial of $y$ of homogeneous degree $d$.
\begin{Definition} {\rm The value of coefficients
$(c_{\alpha_1...\alpha_d})$ is called {\it discriminantal} if for this
value the system of equations
\begin{equation}\label{E1.1}
\frac{\partial{f}}{\partial{y_{\alpha}}}=0,\quad  \alpha=0,...,n
\end{equation}
has a solution in ${\bf P}_n$, i.e. a nonzero one.}
\end{Definition}
If the set of all discriminantal values of $c_{\alpha_1...\alpha_d}$ is an
algebraic manifold of codimension 1 in the space of coefficients, then it is
called {\it the discriminant} of $f$, denoted by $D(f)$. In particular case
when $f$ is a $d-$linear form, its coefficients may be viewed as
elements of a $d-$dimensional matrix $(a_{i_1...i_d})$.

In this paper we show how the invariants of $n-$ary forms can be produced
from the discriminants of multilinear forms (determinants of
multidimensional matricies), which should be considered as the
generalization of the operation of taking classical hessians and resultants.
The algorithm of computation of discriminants of multilinear forms is
considered in paper [1].

\noindent
{\bf Acknoledgment} I want to thank Professor I.M.Gelfand for numerous
critical discussions, which, in particular, helped me to turn aside from
the "naive" definition of a "skew defferential" of a "skew form"
$F=\sum_{i_1<\dots<i_n}f_{i_1\dots i_n}dx_{i_1}\dots dx_{i_n}$ as
$$d_{\varepsilon}F=\sum_{i_1<\dots<i_k<\dots<
i_{n+1}}\varepsilon^{i_k}\frac{\partial f_{i_1\dots{\widehat i_k}\dots
i_{n+1}}}{\partial x_{i_k}}dx_{i_1}\dots dx_{i_{n+1}}$$


\section{Hyperpolarisation and Hyperhessians}
\subsection{Basic example. Binary forms}

{\bf Example.} Let $f=c_{30}x^3+c_{21}x^2y+c_{12}xy^2+c_{03}y^3$ be a
homogeneouse polinomial of 3-rd degree. Take its full polarisation form,
which is a 3-linear form $(a_{i_1i_2i_3})$ with coefficients
$$
a_{111}=\frac{\partial^3{f}}{\partial{x^3}}%
$$
$$
a_{112}=a_{121}=a_{211}=\frac{\partial^3{f}}{\partial{x^2}\partial{y}}%
$$
$$
a_{122}=a_{212}=a_{221}=\frac{\partial^3{f}}{\partial{x}\partial^2{y}}%
$$
$$
a_{222}=\frac{\partial^3{f}}{\partial{y^3}}%
$$
\begin{Proposition}
The discriminant of the full polarisation form $(a_{i_1i_2i_3})$ is equal
to the discriminant of $f$.
\end{Proposition}

\noindent
Now let $f=c_{k0}x^k+c_{{k-1},1}x^{k-1}y^1+...+c_{0k}y^k$ be a homogeneous
polynomial of degree $k$. Let $(a_{i_1...i_k})$ be a $k-$linear form (${%
\underbrace{2\times{...}\times{2}}_{\mbox{\scriptsize{k times}}}}$ matrix)
with coefficients
$$
a_{i_1...i_k}=\frac{\partial^{i+j}f}{\partial{x^i}\partial{y^j}}%
$$
where $j$ is the number of $i_k$ which are equal to 1, $j$ is the number of $%
i_k$ which are equal to 2.

\begin{Proposition}
The discriminant of the full polarisation form $(a_{i_1...i_k})$ is a
product of $GL_2$ invariants of binary form $f$, and in particular is
divisible by the discriminant of $f$.
\end{Proposition} 

\subsection{General Setting}

Throughout the paper let $V$ be a linear vector space of dimension $n+1$
and let ${\bf P}_{n}:={\bf P}(V)$ denote its projectivisation. Let
$f(x_0,x_1,\dots,x_n)$ be a homogeneous polynomial of degree $k$, or a
global section of ${\cal O}_{{\bf P}_{n}}(k)$.  For a sequence of
positive integers ${K}:=(k_1,...,k_d)$ we have a $d$-linear form with
coefficients

\begin{equation} \label{E1.2}
a_{{I}_1...{I}_d}=\frac{\partial
^{|{I}_1+...+ {I}_d|}f}{\partial {\bf {x}}^{{I}_1+...+{I}_d}}
\end{equation}

\noindent
where ${I}_m\in
\{(i_0,...,i_n)|\quad i_0+...+i_n=k_m \}$ are multiidecies. For a
multiindex ${I}=(i_0,...,i_n)$ we denote $\partial {X^I}:= \partial
{x_0^{i_0}}...\partial {x_n^{i_n}}$. Note that $|{I}_1+...+{I}_d|
=k_1+...+k_d$.  Each $a_{{I}_1...{I}_d}$ is again a
homogeneous polynomial of $(x_0,...,x_n)$ but of degree $k-(k_1+...+k_d)$.

\begin{Definition}
{\rm We will call the form $(a_{{I}_1...{I}_d})$
defined by (1.2) {\it ${K}$-polarisation form} (or just a {\it
hyperpolarisation})and denote it as $P^{(k_1...k_d)}(f)$.}
\end{Definition}

So ${K}$-polarisation is a map from the space ${\cal O}(k)$ of
homogeneous polynomials of degree $k$ to the space of $d-$linear forms
with coefficients in ${\cal O}(k-|{K}|)$ which we will
denote by ${\cal O}(k-|{K}|)\otimes {S^{K}T^*(V)}$ (where we use the
notation $S^{K}U :=S^{k_1}U\otimes\dots\otimes S^{k_d}U$):  $${\cal
O}(k)\to {\cal O}(k-|{K}|)\otimes {S^{K}T^*(V)}$$ $$f\mapsto P^{(k_1\dots
k_d)}(f)$$

\noindent
{\bf Example} For $k_1=\dots =k_d=1$ the corresponding hyperpolarisation
form is the form of usual $d$-th polarisation of $f$.

\begin{Definition}
{\rm The discriminant of the form $(a_{{I}_1...{I}_d})=P^{(k_1\dots
k_d)}(f)$ is called $(k_1,...,k_d)-${\it hessian} (or just {\it
hyperhessian}) of $f$ and is denoted as}
$$
{\cal H}^{(k_1...k_d)}(f):=D(a_{{I}_1...{I}_d })=
D(P^{(k_1...k_d)}(f))$$
\end{Definition}

\noindent
{\bf Example}
The usual hessian corresponds to the case when $d=2$ and $k_1=k_2=1$.

\begin{Proposition} Let $k_1+...+k_d=k$. Then the ${K}$-hessian
${\cal H}^{(k_1...k_d)}(f)=D(a_{{\bf  I}_1...{I}_d})$ is a product of
$GL_n$ invariants of $n$-ary form $f$.
\end{Proposition}

So, provided we know how to compute discriminants of $d-$linear forms
(see paper [1] for an outline of the algorithm), each partition
$k_1,...,k_d$ of an integer number $k=k_1+...+k_d$ gives us a set of
invariants of $n-$ary form $f$ of degree $k$. In Section 4 we will give
an example of use of hyperhessians. 

\section{Hyperjacobians and Hyperresultants}

\subsection{Basic Example}

Let $f_1=c^{\prime}_{20}x^2+c^{\prime}_{11}xy+c^{\prime}_{02}y^2$ and
$f_2=c''_{20}x^2+c''_{11}xy+c''_{02}y^2$ be a pair of polynomials.
Take a 3-linear form $(a_{i_1i_2j})$ (with $2\times 2\times 2$ matrix of
coefficients), such that

$$
a_{11j}=\frac{\partial^2{f_j}}{\partial{x^2}}
$$
$$
a_{12j}=a_{21j}=\frac{\partial^2{f_j}}{\partial{x}\partial{y}} $$
$$
a_{22j}=\frac{\partial^2{f_j}}{\partial{y^2}}
$$
i.e. $a_{i_1i_2j}=(P^{(1,1)}(f_j))_{i_1i_2}$.
\begin{Proposition}
The discriminant of the form $(a_{i_1i_2j})$ is equal to the resultant of
$f_1$ and $f_2$ $$D(a_{i_1i_2j})=Res(f_1,f_2)$$
\end{Proposition}

Now let $f_1=c'_{k0}x^k+c'_{{k-1},1}x^{k-1}y^1+...+c'_{0k}y^k$ and
$f_2=c''_{k0}x^k+c''_{{k-1},1}x^{k-1}y^1+...+c''_{0k}y^k$ be a
pair of homogeneous polynomials of degree $k$. Let $(a_{i_1...i_kj})$ be a
$k+1$-linear form with
${\underbrace{2\times{...}\times{2}}_{\mbox{\scriptsize{k+1 times}}}}$
matrix of coefficients
$$
a_{i_1...i_kj}=
\frac{\partial^{j_1+j_2}f_j}{\partial{x^{j_1}}\partial{y^{j_2}}}
$$
where $j_1$ is the number of $i_k$ which are equal to 1, $j_2$ is the
number of $i_k$ which are equal to 2.

\begin{Proposition} The
discriminant of the form $(a_{i_1...i_kj})$ is a product
of $GL_2$ invariants of a pair of binary forms $f_1,f_2$, and in
particular is divisible by the resultant $Res(f_1,f_2)$.
\end{Proposition}

\subsection{General Setting}

Let ${M}:=(m_1,...,m_{d_1})$ be a sequence of positive integers and
let $F^{M}(k):=\{f_{I}\}_{{I}\le{M}}$ be a system of
homogeneous polynomials of degree $k$ ennumerated by
the indecies ${I}=(i_1,\dots,i_{d_1})$, or a global section of
${\cal O}(k)^{M}:={\cal O}(k)^{\oplus m_1}\otimes\dots\otimes {\cal
O}(k)^{\oplus m_{d_1}}$.

Let ${K}:=(k_1,...,k_{d_2})$. For each $f_{I}$ take its $K$-polarisation
form $P^{(k_1...k_d)}(f_{I})$ as described in Section 1.2. Then we get a
$(d_1+d_2)$-linear form $(a_{i_1\dots i_{d_1}{J}_1...{J}_{d_2}})$ such
that

\begin{equation} \label{E2.3}
a_{i_1\dots i_{d_1}{J}_1...{J}_{d_2}}=
(P^{(k_1...k_d)}(f_{i_1\dots i_{d_1}}))_{{J}_1\dots{J}_{d_2}}
\end{equation}
where as in Section 1.2 ${J}_i\in
\{(j_0,...,j_n)|\quad j_0+...+j_n=k_i\}$ are multiidecies. Each
$a_{i_1\dots i_{d_1}{J}_1...{J}_{d_2}}$ is a homogeneous polynomial of
$(x_0,\dots,x_n)$ of degree $k-|K|$.

\begin{Definition}
{\rm The form defined by (2.3) is called $K$-{\it Jacobi form}  of the
section $(f_{i_1\dots i_{d_1}})_{I\le M}$, considered as a map of ${\bf
P}_n$ to the space of $n_1\times\dots\times n_{d_1}$ tensors (or a tensor
field on ${\bf P}_n$) and is denoted by}
$$J^K(F^M)=\frac{{\rm D}^K((f_I)_{I\le M})}{{\rm
D}(x_1,\dots,x_n)^{|K|}}$$

{\rm The discriminant of $K$-Jacobi form is called $K$-{\it Jacobian} of
$(f_{i_1\dots i_{d_1}})_{I\le M}$ (or just a {\it hyperjacobian}).}

\end{Definition}

So taking the $K$-Jacobi form of tensor fields $F^M\in{\cal O}(k)^{M}$
gives us a map
$$J^K:{\cal O}(k)^{M}\to{\cal O}(k-|K|)^{M}\otimes S^{K}T^*(V)$$

\noindent
{\bf Example} For $F^M=(f_1,\dots,f_m)$ and $K=(1)$ the
hyperpolarisation form $P^{(1)}(f_i)$ is the differential $df_i$ and the
corresponding $K$-Jacobi form is the usual Jacobi matrix of the map ${\bf
P}_n\to{\bf P}_m$.

\subsection{Hyperresultant}

Let $F^M$ be a section of ${\cal O}(k)^{\oplus m}$, i.e. is a map ${\bf
P}_n\to {\bf P}_m$ of degree $k$, or a system $(f_1,\dots,f_m)$ of $m$
homogeneous polynomials of degree $k$. For $K=\underbrace{(1,\dots,1)}_{k
\ {\rm times}}$ the hyperpolarisation form $P^K(f_i)$ is the form with
constant coefficients.  For the corresponding $K$-Jacobian we will reserve
a special name.

\begin{Definition} {\rm The $\underbrace{(1,\dots,1)}_{k\ {\rm
times}}$-Jacobian of a system $(f_1,\dots,f_m)$ of polynomials of degree
$k$ is called {\it hyperresultant} of this system and is denoted by
$$R_m(f_1,...,f_m):=D((a_{i_1\dots i_kj})_
{
i_1,\dots,i_k\le n
\  j\le m
})=D(P^{(1,\dots,1)}(f_j)_{j\le m})$$}.

\end{Definition}

\begin{Proposition}
For $m=2$ the hyperresultant $R_2(f_1,f_2)$ is divisible by the usual
resultant $R(f_1,f_2)$.

\end{Proposition}

\noindent
{\bf Example} Let $f_1=a_{20}x^2+a_{11}xy+a_{02}y^2,
\  f_2= b_{20}x^2+b_{11}xy+b_{02}y^2$ and $f_3=
c_{20}x^2+c_{11}xy+c_{02}y^2$. If we write $f_i(x,y)$ as functions of a
nonhomogeneous variable $x^{\prime}:= \frac{x}{y}$ then the Wronskian $$
W(f_1,f_2,f_3)=\left|{
\matrix{
\frac{\partial^2{f_1}}{\partial{x'^2}} &
\frac{\partial{f_1}}{\partial{x'}} & f_1 \cr
\frac{\partial^2{f_2}}{\partial{x'^2}} &
\frac{\partial{f_2}}{\partial{x'}} & f_2 \cr
\frac{\partial^2{f_3}}{\partial{x'^2}} &
\frac{\partial{f_3}}{\partial{x'}} & f_3 \cr }}\right| $$
is equal to 0 iff $f_1,f_2$ and $f_3$ are linearly dependent. On the
other hand for each $f_i$ we take the $2\times 2$ matrix $P^{(1,1)}(f_i)$
and make from these matricies a $2\times 2\times 3$ form
$R_3(f_1,f_2,f_3)$ according to (2.3).
Then for 3-resultant of $f_1, f_2$ and $f_3$ (the discriminant of
$J^{(1,1)}(f_1,f_2,f_3)$) there is an equality
$$R_3(f_1,f_2,f_3)=W(f_1,f_2,f_3)^2$$

\subsection{Jacobi Sequence}
Let $K=(1)$. Then the corresponding 1-Jacobi map defines a sequence
$${\cal O}(k)\textmap{J^{(1)}}{\cal O}(k-1)\otimes
T^*(V)\textmap{J^{(1)}}\dots \textmap{J^{(1)}}{\cal O}(0)\otimes
T^*(V)^{\otimes k}$$

Now let us extend the Jacobi map from homogeneous polynomials
to the space ${\cal O}$ of analytic (or for some purposes just
differentiable) functions on an $n$-dimensional local chart $X$ and denote
$T^*(X)\cong{\cal O}\otimes T^*_x(X)$ (for any given $x\in X$) the space
of analytic sections of cotangent bundle of $X$.  This gives us an
infinite sequence:  $${\cal
O}\textmap{J^{(1)}}T^*(X)\textmap{J^{(1)}}\dots
\textmap{J^{(1)}}T^*(X)^{\otimes k}\textmap{J^{(1)}}\dots$$
which we will call a {\it Jacobi sequence}.

\section{Gramm Complexes}

\subsection{Combinatorial Remark}

Let $d>0$ be an integer and $k_1\ge\dots\ge k_p$ be its partition
$k_1+\dots+k_p=d$. The group $S_d/S_{k1}\times\dots\times
S_{k_p}$ has a free representation on the set of words of length $d$ of
$p$ letters. There is a noncanonical but, nevertheless, natural
lexicographic ordering on this representation space, which induces an
oredering on the group. For $\sigma\in S_d/S_{k1}\times\dots\times
S_{k_p}$ let ${\rm ord}(\sigma)$ denote the ordinal number of the group
element $\sigma$ with respect to this ordering.

For an integer $d$ let $\varepsilon$ be the
primitive $d!$-th root of unity $\varepsilon^{d!}=1$.

\noindent
{\bf Notation } For an $n$-dimensional vector space $V$ let
$V^{\otimes d}_{k}\subset V^{\otimes d}$ denote the subspace
of $d$-vectors with the condition of $\varepsilon^k$-{\it skew symmetry}
of their coordinates:  $$a_{\sigma(i_1\dots i_d)}=\varepsilon^{k\
\scriptsize{ord}(\sigma)}a_{i_1\dots i_d}$$

Here each set of indices $i_1,\dots,i_d$ is arranged as:
$$i_1=\dots=i_{k_1}$$
$$i_{k_1+1}=\dots=i_{k_1+k_2}$$
$$\dots$$
$$i_{k_1+\dots+k_{d-1}+1}=\dots=i_{d}$$
with $k_1\ge\dots\ge k_d$, and $\varepsilon^{k\frac{d!}{k_1!\dots
k_d!}}=1$. Then such a set may be considered as the initial (in
lexicographic ordering) element of representation space of the group
$S_d/S_{k1}\times\dots\times S_{k_p}\ni\sigma$ with the corresponding
well defined ordinal function ${\rm ord}(\sigma)$ on it.

\begin{Proposition}
$$V^{\otimes d}=\bigoplus_{k=1}^{k=d!}V^{\otimes d}_{k}$$

\end{Proposition}

\noindent
{\bf Example} Let $\dim V=3$, $d=3$. Then
$$V^{\otimes 3}=V^{\otimes 3}_0\oplus V^{\otimes 3}_1\oplus V^{\otimes
3}_5\oplus V^{\otimes 3}_2\oplus V^{\otimes 3}_4\oplus V^{\otimes 3}_3.$$
Here $V^{\otimes 3}_0$ is symmetric part of dimension 10 with nonzero
components with idices sets $(111),(222),(333),$
$(112),(113),(223),(221),(331),(332),(123)$, $V^{\otimes 3}_1$ and
$V^{\otimes 3}_5$ are parts of dimension 1 with index set $(123)$,
$V^{\otimes 3}_2$ and $V^{\otimes 3}_4$ are parts of dimension $7$ with
index set $(112),(113),(223),(221),(331),(332),(123)$, $V^{\otimes 3}_3$
is antisymmetric part of dimension $1$ with index set $(123)$. An example
of the generic element of $V^{\otimes 3}_2$ is
$$
a_{112}e_1\otimes e_1\otimes e_2+
\varepsilon^2a_{112}e_1\otimes e_2\otimes e_1+
\varepsilon^4a_{112}e_2\otimes e_1\otimes e_1+$$
$$
a_{113}e_1\otimes e_1\otimes e_3+
\varepsilon^2a_{113}e_1\otimes e_3\otimes e_1+
\varepsilon^4a_{113}e_3\otimes e_1\otimes e_1+$$
$$
a_{223}e_2\otimes e_2\otimes e_3+
\varepsilon^2a_{223}e_2\otimes e_3\otimes e_2+
\varepsilon^4a_{223}e_3\otimes e_2\otimes e_2+$$
$$
a_{221}e_2\otimes e_2\otimes e_1+
\varepsilon^2a_{221}e_2\otimes e_1\otimes e_2+
\varepsilon^4a_{221}e_1\otimes e_2\otimes e_2+$$
$$
a_{331}e_3\otimes e_3\otimes e_1+
\varepsilon^2a_{331}e_3\otimes e_1\otimes e_3+
\varepsilon^4a_{331}e_1\otimes e_3\otimes e_3+$$
$$
a_{332}e_3\otimes e_3\otimes e_2+
\varepsilon^2a_{332}e_3\otimes e_2\otimes e_3+
\varepsilon^4a_{332}e_2\otimes e_3\otimes e_3+$$
$$
a_{123}e_1\otimes e_2\otimes e_3+
\varepsilon^2a_{123}e_1\otimes e_3\otimes e_2+
\varepsilon^4a_{123}e_2\otimes e_1\otimes e_3+
\varepsilon^0a_{123}e_2\otimes e_3\otimes e_1+
\varepsilon^2a_{123}e_3\otimes e_1\otimes e_2+
\varepsilon^4a_{123}e_3\otimes e_2\otimes e_1
$$
where $\varepsilon$ is the primitive $3!$-th root of unity.

{\bf Notation} For a given basis in $V$ denote by $p_k$ the projection map
of $V^{\otimes d}$ onto $\varepsilon^k$-skew symmetric component
$V^{\otimes d}_k$ for $k\le d!$: $$p_k:V^{\otimes d}\to V^{\otimes d}_k.$$

\subsection{Gramm Forms}
Let $F\in U^{*\otimes d}$ be a $d$-linear form. Then for
each $\{u_1,\dots,u_m\}\in\underbrace{U\times\dots\times U}_{d\ {\rm
times}}$ we have a set of expressions
$$<u_{i_1},\dots,u_{i_d}>:=F(u_{i_1},\dots,u_{i_d}),\quad
i_1,\dots,i_d=1,...,m.$$ Let $D(<u>)$ be the discriminant of
the $\underbrace{m\times\dots\times m}_{d\ {\rm times}}$ form with
coefficients $<u_{i_1},\dots,u_{i_d}>$.  Then $D(<u>)^{m/d{\rm
deg}D(<u>)}$ is an expression of $u_1,...,u_m$ which has homogeneous
degree $1$ with respect to each of $u_i$.

\begin{Definition}

{\rm The function $$G^d:\{u_1,\dots,u_m\}\mapsto D(<u>)^{\frac{m}{d{\rm
deg}D(<u>)}}$$ is called the $d$-{\it Gramm form}.}

\end{Definition}

\noindent
{\bf Example} Let $d=2$. Then $(<u_{i_1},u_{i_2}>)_{1\le i_1,i_2\le m}$ is
the usual Gramm matrix of the $m$-tuple of vectors $u_1,\dots,u_m$.

\begin{Proposition}
Let $g\in GL(m)$. Then $G^d(<g(u)>)=\mid g\mid G^d(<u>)$.
\end{Proposition}
\noindent
{\bf Proof} This is the consequence of the fact that the discriminant of
$\underbrace{m\times\dots\times m}_{d\ {\rm times}}$ form is a $GL(m)$
invariant.

\begin{Proposition}
If $u_1,\dots,u_m$ are linearly dependent then $G^d(<u>)\equiv 0$.
\end{Proposition}

In particular, if the number of vectors $u_i$ is less then the
"dimensionality" of the form $m<d$ then $G^d(u_1,\dots,u_m)\equiv 0$.

Let $X$ be a manifold, and let $F_x$ be a section of $T^*(X)^{\otimes
d}$,  i.e.  a field of $d-$linear forms on $T(X)^{\times m}$. Then the map
$G^d(F_x):u_1,\dots,u_m\mapsto
G^d(u_1,\dots,u_m;F_x)$ gives us a measure of integration on
any $d$-submanifold of $X$.  Let $\mu(X)$ denote the space of analytic
measures on $X$. Thus taking a Gramm form defines a map:
$$T^*(X)^{\otimes d}\to \mu(X)$$ $$F(\bullet)\mapsto G^d(\bullet;F)$$

\subsection{Gramm Complex}

According to Section 3.1 the  $\underbrace{m\times\dots\times m}_{d\ {\rm
times}}$ form $(<u_{i_1},\dots,u_{i_d}>)$ may be decomposed into
$\varepsilon^k$-skew symmetric parts $(<u_{i_1},\dots,u_{i_d}>)_k$:

$$\begin{array}{ccc}
u_1\times\dots\times u_m&\textmap{}&(<u_{i_1,\dots,u_{i_d}}>)\\
&&\downarrow^{p_k}\\
&&(<u_{i_1,\dots,u_{i_d}}>)_k
\end{array}\quad {\rm for}\  k=1,...,d!$$

\begin{Definition}
Functions $G^d_{k}: u_1,\dots,u_m\mapsto
D((F(u_{i_1},\dots,u_{i_d}))_k)$ are called $k$-{\it skew Gramm forms}.
\end{Definition}

If $X$ is a manifold and $F$ is a field of
$d-$linear forms, i.e. a section of $T^*(X)^{\otimes d}$, then each
$G^d_k,\quad k=1,\dots,d!$ gives us a measure on $X$. Thus
taking a $k$-skew Gramm forms defines a set of maps:
$$\begin{array}{ccc}
T^*(X)^{\otimes d}&\textmap{G^d_1}&\mu(X)\\
&\textmap{G^d_2}&\mu(X)\\
&\dots&\\
&\textmap{G^d_{d!}}&\mu(X)\\
\end{array}$$

Each image $G^d_k(T^*(X)^{\otimes d})$ is a linear subspace $\mu^d_k(X)$
in the space $\mu(X)$.

For a fixed $k$ the map $G^{\bullet}_k$ applied to the Jacobi sequence:
$$\begin{array}{ccccccccccc}
{\cal O}&\textmap{J^{(1)}}&T^*(X)&\textmap{J^{(1)}}&\dots&
\textmap{J^{(1)}}&T^*(X)^{\otimes k}&\textmap{J^{(1)}}& T^*(X)^{\otimes
k+1}&\textmap{J^{(1)}}&\dots\\
\downarrow^{G^0_k}&&\downarrow^{G^1_k}&&\dots&&\downarrow^{G^k_k}&&
\downarrow^{G^{k+1}_k}&&\dots\\
0&&\dots&&0&&\mu^k_k(X)&&\mu^{k+1}_k(X)&\dots\\

\end{array}$$
gives us an infinite sequence:
$$\begin{array}{ccccccccccccccc}
0&\textmap{\delta_k}&\dots&\textmap{\delta_k}&0&
\textmap{\delta_k}&\mu^k_k(X)&\textmap{\delta_k}&\mu^{k+1}_k(X)
&\textmap{\delta_k}&\dots&\textmap{\delta_k}
&\mu^n_k(X)&\textmap{\delta_k}\dots

\end{array}$$

\begin{Proposition}
For integer $1\le k< d!$ the operator $\delta^k$ has a final order:
$$\delta^d\equiv 0.$$

\end{Proposition}

\begin{Definition}
We call the operator $\delta_k$ $k-${\it differential} and the
sequence () $k-${\it Gramm complex}.
\end{Definition}

\noindent
{\bf Examples}

\noindent
1. For $k=d!$ (or, equivalently, $k=0$) the 0-defferential corresponds to
symmetric part of Gramm form and $\delta_k$ has infinite order.

\noindent
2. For $k=d!/2$ the $k$-Gramm complex is just the de Rham complex of our
manifold $X$.

So the set of $k$-Gramm complexes may be considered as a deformation
connecting de Rham complex with differential $\delta^2\equiv 0$ with the
complex (Example 1.) with differential of infinite order.


\section{Applications of Hyperjacobians}

Let $f=c_{40}x^4+a_{31}x^3y+c_{22}x^2y^2+c_{13}x^1y^3+c_{04}y^4$ be a
binary form of degree 4. It is known (see [2])
that the ring of $GL_2$ invariants of
$f$ is generated by the polynomials: 1) Hankel determinant
$$
{\cal {H}}:={\frac{1}{8}}\left|{
\matrix{
24c_{40} & 6c_{31} & 4c_{22}  \cr
6c_{31}   & 4c_{22} & 6c_{13}  \cr
4c_{22}   & 6c_{13} & 24c_{04}
}}\right|
$$
2) apolara
$$
{\cal {A}}:=c_{22}^2+3c_{31}c_{13}+12c_{40}c_{04}
$$
Denote also by ${\cal {D}}$ another invariant from this ring - the
discriminant of $f$.\\ On the other hand we can take a (1,1,1)-hessian of
$f$, denoted by $f^{111}:=H_{111}(f)$. $f^{111}$ is again a polynomial in
$(x,y)$ of degree 4.
\begin{Proposition}
$$D(f^{111})=2^{36}3^6{\cal{D}}{\cal{H}}^{6}$$
$$R(f,f^{111})=2^{24}3^{12}{\cal{D}}^2{\cal{A}}^4$$
\end{Proposition}
So as soon as we know how to compute the hyperhessian $H_{111}(f)$ (i.e.
the 3-dimensional determinants, see [1] for the outline of an algorithm)
we can get the generators of the ring of invariants by taking the
classical discriminants and resultant of $f$ and $f^{111}$.

The experiments with polynomials and their systems of higher degree and of
more variables are left to the reader.

\end{document}